\documentclass[12pt,reqno]{amsart}

\usepackage[latin1]{inputenc}
\usepackage[T1]{fontenc}
\usepackage[british]{babel}
\usepackage{amssymb,graphicx,url}
\usepackage[margin=30pt,font=small]{caption}
\usepackage[margin=0pt]{subfig}

\theoremstyle{plain}
\newtheorem{theorem}{Theorem}[section]
\newtheorem{proposition}[theorem]{Proposition}
\newtheorem{lemma}[theorem]{Lemma}
\newtheorem{corollary}[theorem]{Corollary}

\theoremstyle{definition}
\newtheorem{definition}[theorem]{Definition}

\theoremstyle{remark}
\newtheorem{example}[theorem]{Example}
\newtheorem{remark}[theorem]{Remark}
\newtheorem*{acknowledgements}{Acknowledgements}

\numberwithin{equation}{section}
\numberwithin{figure}{section}

\newcommand{\Lie}[1]{\operatorname{\textsl{#1}}}
\newcommand{\lie}[1]{\operatorname{\mathfrak{#1}}}
\newcommand{\GL}{\Lie{GL}}
\newcommand{\Sl}{\lie{sl}}

\makeatletter
\newcommand{\abs}{\@ifstar{\absManualScale}{\absAutoScale}}
\newcommand{\absAutoScale}[1]{\left\vert#1\right\vert}
\newcommand{\absManualScale}[2][]{\mathopen #1\lvert #2 \mathclose #1\rvert}

\newcommand{\inp}{\@ifstar{\inpManualScale}{\inpAutoScale}}
\newcommand{\inpAutoScale}[2]{\left< #1, #2 \right>}
\newcommand{\inpManualScale}[3][]{\mathopen #1\langle #2, #3 
\mathclose #1\rangle} 
\makeatother

\DeclareMathOperator{\re}{Re}
\DeclareMathOperator{\diag}{diag}
\DeclareMathOperator{\image}{im}
\DeclareMathOperator{\rank}{rank}
\DeclareMathOperator{\Span}{span}
\DeclareMathOperator{\Stab}{Stab}
\DeclareMathOperator{\CInt}{CInt}
\DeclareMathOperator{\Int}{Int}
\DeclareMathOperator{\Aut}{Aut}

\newcommand{\ii}{\mathbf i}

\begin{document}
\title{Toric Hypersymplectic Quotients}

\author{Andrew Dancer}
\address[Dancer]{Jesus College, Oxford, OX1 3DW, United Kingdom}
\email{dancer@maths.ox.ac.uk}

\author{Andrew Swann}
\address[Swann]{Department of Mathematics and Computer Science\\
University of Southern Denmark\\
Campusvej 55\\
DK-5230 Odense M\\
Denmark}
\email{swann@imada.sdu.dk}

\begin{abstract}
  We study the hypersymplectic spaces obtained as quotients of flat
  hypersymplectic space~\( \mathbb R^{4d} \) by the action of a compact
  Abelian group.  These \( 4n \)-dimensional quotients carry a
  multi-Hamilitonian action of an \( n \)-torus.  The image of the
  hypersymplectic moment map for this torus action may be described by a
  configuration of solid cones in~\( \mathbb R^{3n} \).  We give precise
  conditions for smoothness and non-degeneracy of such quotients and show
  how some properties of the quotient geometry and topology are constrained
  by the combinatorics of the cone configurations.  Examples are studied,
  including non-trivial structures on~\( \mathbb R^{4n} \) and metrics on
  complements of hypersurfaces in compact manifolds.
\end{abstract}

\subjclass[2000]{Primary 53C25; Secondary 53D20, 53C55, 57S15}

\keywords{Hypersymplectic structure, neutral hyperKähler manifold, toric
variety, moment map}

\maketitle

\newpage
\tableofcontents

\newpage
\section{Introduction}
\label{sec:introduction}
An important construction in symplectic geometry is the \emph{symplectic
quotient} of Marsden and Weinstein.  Given a symplectic action of a Lie
group \( G \) on a symplectic manifold \( M \), this produces, under fairly
mild hypotheses, a new symplectic manifold of dimension \( \dim M - 2 \dim
G \).

One particularly nice class of symplectic examples is that of toric
varieties.  Delzant~\cite{De} and Guillemin~\cite{Gu}, have shown that a
large class of toric varieties may be produced as symplectic quotients of
\( \mathbb C^d \) with its flat K\"ahler structure by a subtorus \( N \)
of~ \( \mathbb T^d \).  Their work also shows that the symplectic and
K\"ahler geometry of these \( 2n \)-dimensional examples is determined by
concrete descriptions of the moment polytope as an intersection of closed
half-spaces in~\( \mathbb R^n \).

In the present paper we shall investigate an analogous construction in a
situation where there are multiple symplectic structures, namely for
hypersymplectic manifolds.  We consider the geometries arising from a
hypersymplectic quotient construction for compact Abelian subgroups~\( N \)
of \( \mathbb T^d \) acting on \( \mathbb C^{d,d}=\mathbb C^d\times \mathbb
C^d \).  We determine conditions for the quotients to be smooth manifolds
and to admit non-degenerate geometric structures in two ways: firstly by
direct considerations, and then by using the (hypersymplectic) moment map
for the action of a torus on the quotient.  The closed half-spaces of
Delzant and Guillemin are now replaced by solid cones in \( \mathbb R^{3n}
\) and we demonstrate how properties of the quotient may be deduced from
particular descriptions of such cone configurations.

\emph{Hypersymplectic structures} were defined in a paper of
Hitchin~\cite{H} and have appeared in recent works such
as~\cite{Hu,Ka,FPPW,AD}.  One has an \emph{indefinite} K\"ahler metric \( g
\) of signature \( (2n, 2n) \), together with a covariant constant
endomorphism \( S \) of the tangent bundle, such that \( S^2 \) equals the
identity, \( S \) anti-commutes with the complex structure \( I \), and \(
g(SX, SY)=-g(X,Y) \).  Now \( I, S \) and \( T = IS \) generate an action
of the Lie algebra \( \Sl(2, \mathbb R) \) on each tangent space.
Moreover, \( I,S,T \), together with the metric \( g \), define three
symplectic forms by
\begin{gather*}
  \omega_I(X, Y) = g(X, IY), \qquad
  \omega_S(X,Y)  = g(X, SY), \\
  \omega_T(X, Y) = g(X, TY),
\end{gather*}
hence the name ``hypersymplectic''.  Every hypersymplectic manifold is
\emph{neutral Calabi-Yau}, that is, Ricci-flat K\"ahler with signature \(
(2n,2n) \).  Hypersymplectic manifolds are split-quaternion analogues of
hyperk\"ahler manifolds and are sometimes referred to as ``neutral
hyperk\"ahler manifolds''.

Hitchin described a quotient construction for hypersymplectic manifolds
in~\cite{H}, analogous to the hyperk\"ahler quotient of~\cite{HKLR}.  If \(
G \) acts preserving a hypersymplectic structure, then under mild
conditions, we have, for each \( X \in \lie g \), a moment map \( \mu^X \)
taking values in \( \mathbb R^3 \).  This map satisfies the defining
equation
\begin{equation}
  \label{eq:mm}
  d \mu^X (Y) = (\omega_I (X, Y), \omega_S (X, Y), \omega_T (X,Y)),
\end{equation}
where we identify \( X \) with the vector field it induces via the group
action.  Of course, equation~\eqref{eq:mm} only gives~\( \mu^X \) up to an
additive constant; these constants are partially restricted by the
additional assumption that the maps \( \mu^X \) combine to define a \( G
\)-equivariant map \( \mu \) taking values in \( \lie g^* \otimes \mathbb
R^3 \).  If \( G \)~is Abelian, then any choice of \( \mu^X \) gives an
equivariant map.  The hypersymplectic quotient is now defined to be \(
\mu^{-1}(0)/G \).  When \( \mu \) has maximal rank and the action of~\( G
\) is free, the quotient has dimension \( 4\dim G \) less than the original
hypersymplectic manifold.  It inherits closed two-forms from \( \omega_I
\), \( \omega_S \) and \( \omega_T \) and one expects these to define a
hypersymplectic structure; however degeneracies may occur on a certain
locus in the quotient.

In this paper we shall concentrate on hypersymplectic quotients of flat
space \( \mathbb C^{d,d} \) by compact Abelian groups, although we prove
some general results controlling smoothness and non-degeneracy of arbitrary
hypersymplectic quotients.  In some ways, the picture is intermediate
between that of K\"ahler quotients and of hyperk\"ahler quotients as
studied in~\cite{BD}.  In the hyperk\"ahler case the quotients are
necessarily non-compact, whereas in the K\"ahler case many quotients are
compact.  For the hypersymplectic situation, we show how to produce
non-compact non-singular structures on \( \mathbb R^{4n} \) that are not
flat.  It is also easy to produce compact quotient sets and with a little
more work smooth examples may be found, but on the other hand these always
have singularities of the hypersymplectic structure.  We show that these
compact quotients produce non-degenerate structures on hypersurface
complements in real analytic subvarieties of compact toric varieties.  All
the quotients we produce carry a natural involution.  We discuss in detail
some particular examples, such as the hypersymplectic analogues of the
Calabi and Gibbons-Hawking multi-instanton spaces, and all examples
obtainable as quotients of \( \mathbb C^{2,2} \) by a one-dimensional
group.

\begin{acknowledgements}
  Both authors are members of the \textsc{Edge}, Research Training Network
  \textsc{hprn-ct-\oldstylenums{2000}-\oldstylenums{00101}}, supported by
  The European Human Potential Programme.
\end{acknowledgements}

\section{The flat hypersymplectic structure}
\label{sec:flat}
Our examples will be hypersymplectic quotients of the following flat
hypersymplectic structure.  Let \( \mathbb C^{d,d} \) be \( \mathbb C^d
\times \mathbb C^d \) with the complex structure
\begin{equation*}
  I(z,w) =(z \ii,-w \ii),
\end{equation*}
where \( \mathbf i = \sqrt{-1} \), and with the indefinite K\"ahler metric
\begin{equation*}
  g =\re\Bigl(\sum_{k=1}^d dz_k d \bar z_k - dw_k d \bar w_k \Bigr).
\end{equation*}
Note that \( I \) is not the standard complex structure \( I_0 \) on this
space, which is instead defined by \( I_0 \colon (z,w) \mapsto (z \ii, w
\ii) \).  To distinguish \( I \) from \( I_0 \), we shall refer to this
space with complex structure \( I_0 \) as \( \mathbb C^{2d} \).

If we define \( S(z,w)=(w,z) \), then
\begin{gather*}
  IS(z,w) = I(w,z) = (w \ii,-z \ii),\\
  SI(z,w) = S(z \ii,-w \ii) = (-w \ii,z \ii),
\end{gather*}
so \( IS=-SI \).  We define \( T=IS \), so that
\begin{equation*}
  I^2 = -1,\quad S^2=T^2=1,\quad IS=T=-SI.
\end{equation*}

We have the following symplectic forms:
\begin{align*}
  \omega_I &= \frac1{2\ii}\sum_{k=1}^d \left(dz_k\wedge d\bar z_k +
    dw_k\wedge d\bar w_k\right),\\
  \omega_S &= \frac12 \sum_{k=1}^d \left(dz_k \wedge d\bar w_k - dw_k\wedge
    d\bar z_k \right), \\
  \omega_T &= \frac1{2\ii}\sum_{k=1}^d \left(dz_k \wedge d\bar w_k +
    dw_k\wedge d\bar z_k \right).
\end{align*}
Note that \( \omega_S + \ii \omega_T = \sum_{k=1}^d dz_k\wedge d\bar
w_k \), which is a holomorphic \( (2,0) \)-form with respect to \( I \),
but is of type \( (1,1) \) for~\( I_0 \).

\section{Moment maps}
\label{sec:moment}
The torus \( \mathbb T^d \) acts on \( \mathbb C^{d,d} \) by
\begin{equation*}
  (z_k,w_k) \longmapsto (e^{\ii \theta_k}z_k,e^{\ii \theta_k}w_k).
\end{equation*}
This action commutes with \( I \) and \( S \) and hence with \( T \), and
preserves~\( g \).

The moment maps from equation~\eqref{eq:mm} are
\begin{gather*}
  \mu_I \colon (z,w) \longmapsto  \sum_{k=1}^d
  \tfrac12(\abs{z_k}^2+\abs{w_k}^2)e_k + \tilde c_1,\\
  \mu_S + \ii \mu_T \colon (z,w)  \longmapsto \sum_{k=1}^d
  \ii z_k\bar w_k e_k + \tilde c_2 + \ii \tilde c_3,
\end{gather*}
where \( d\mu_I^X(Y)=\omega_I(X,Y) \), etc., \(e_1, \dots, e_d \) are the
standard basis vectors for \(\mathbb R^d \), and \( \tilde c_1,
\tilde c_2, \tilde c_3 \) are arbitrary constant vectors in \( \mathbb
R^d \).

The form of this differs from the hyperk\"ahler moment map~\cite{BD} for
the \( \mathbb T^d \) action on~\(\mathbb H^d \), in that we have \(+
\abs{w_k}^2 \) rather than \(-\abs{w_k}^2 \) in the formula for \( \mu_I
\).  We also recall, for comparison, that the K\"ahler moment map for the
action of \(\mathbb T^d \) on \(\mathbb C^d \) is \( \mu \colon z \mapsto
\sum_{k=1}^d \frac12\abs{z_k}^2 e_k + c\).

As in~\cite{Gu,BD} one considers a compact Abelian subgroup~\( N \) of~\(
\mathbb T^d \) with Lie algebra~\( \lie n \).  We shall take \( \lie n \)
to be the kernel of a surjective linear map~\( \beta\colon \mathbb R^d\to
\mathbb R^n \) given by
\begin{equation*}
  \beta \colon e_k \longmapsto u_k,
\end{equation*}
with \( u_i\in \mathbb Z^n \).  In particular, \( \mathbb R^n \)~is spanned
by \( u_1, \dots, u_d \).  Then we have an exact sequence
\begin{equation}
  \label{eq:beta}
  0 \longrightarrow  {\lie n} \overset\iota\longrightarrow \mathbb R^d
  \overset\beta\longrightarrow \mathbb R^n \longrightarrow 0
\end{equation}
and \( N \)~is defined to be the kernel of the map \(
\exp\circ\beta\circ\exp^{-1}\colon \mathbb T^d \to \mathbb T^n \).  (The
requirement that \( u_i \) be integral exactly guarantees that this
composition is well-defined.)  For a given~\( N \), the map \( \beta \) is
unique up to composition with an element of~\( \Aut(\mathbb Z^n) \), or in
matrix terms up to multiplication by an element of \( \GL(n,\mathbb Z) = \{
A \in M_n(\mathbb Z): \det A\ne0\ \text{and}\ A^{-1}\in M_n(\mathbb Z)\}\).

We let \( \inp\cdot\cdot \) denote the standard inner product with
respect to which \( e_1, \dots, e_d \) are orthonormal.  For each choice
of scalars \( \lambda_1, \dots, \lambda_d \), the set of vectors \(
\{u_1,\dots,u_d\} \) defines a convex polyhedron in \( \mathbb R^n \) by
the equations
\begin{equation}
  \label{eq:poly}
  \inp s{u_k}  \geqslant \lambda_k,  \qquad\text{for \( k=1, \dots, d\).}
\end{equation}
In general, this polyhedron may be non-compact.

There is an exact sequence dual to~\eqref{eq:beta}
\begin{equation*}
  0 \longrightarrow {\mathbb R^n}^* \overset{\beta^*}{\longrightarrow}
  {\mathbb R^d}^* \overset{\iota^*}{\longrightarrow} \lie n^* \longrightarrow
  0.
\end{equation*}
We shall identify \( {\mathbb R^d}^* \) with \( \mathbb R^d \) using \(
\inp\cdot\cdot \).  Now \( \beta^* \) is given by
\begin{equation}
  \label{eq:bstar}
  \beta^* (a) = \sum_{k=1}^d \inp a{u_k} e_k
\end{equation}
and the moment map for \( N \) becomes
\begin{gather*}
  \mu_I \colon (z,w) \longmapsto \sum_{k=1}^d
  \tfrac12(\abs{z_k}^2+\abs{w_k}^2) \alpha_k  + c_1,\\
  \mu_S + \ii \mu_T \colon (z,w) \longmapsto \sum_{k=1}^d \ii z_k\bar w_k
  \alpha_k + c_2+ \ii c_3,
\end{gather*}
where \( \alpha_k = \iota^* e_k \).  We write \(c_j = \sum_{k=1}^d
\lambda_k^{(j)} \alpha_k \) for some scalars \( \lambda_k^{(j)} \).

A point \( (z,w) \) lies in \( \mu_I^{-1}(0) \) if and only if
\begin{equation}
  \label{eq:level}
  \iota^* \Bigl( \sum_{k=1}^d \bigl(\tfrac12(\abs{z_k}^2 +
\abs{w_k}^2)  + \lambda_k^{(1)}\bigr) e_k \Bigr)=0.
\end{equation}
However \(\ker \iota^* = \image \beta^* \), so using
equation~\eqref{eq:bstar}, we see that \eqref{eq:level} is equivalent to
the existence of \( a \in \mathbb R^n \) such that
\begin{equation}
  \label{eq:a}
  \inp a{u_k} = \tfrac12 (\abs{z_k}^2 + \abs{w_k}^2)  +
  \lambda_k^{(1)}, \qquad\text{for \( k=1, \dots, d \).}
\end{equation}
Similarly, \( (z,w ) \in (\mu_S + \ii \mu_T)^{-1}(0) \) if and only if
\begin{equation}
  \label{eq:b}
  \inp b{u_k} = \ii z_k \bar w_k + \lambda_k^{(2)} + \ii
  \lambda_k^{(3)}, \qquad\text{for \( k=1, \dots, d \),}
\end{equation}
for some \( b \in \mathbb C^n \).

Equations \eqref{eq:a} and~\eqref{eq:b} give a description of the level set
\( \mu^{-1}(0) \).  The hypersymplectic quotient \( M \) of \( \mathbb
C^{d,d} \) by~\( N \) is defined to be
\begin{equation*}
  M = \mu^{-1}(0)/N.
\end{equation*}
This is a Hausdorff topological space; as we will see it may or may not be
a smooth manifold.

\section{Non-degeneracy of the quotient geometry}
\label{sec:non-degeneracy}
In this section we shall consider when hypersymplectic quotients are
smooth, and at which points the hypersymplectic structure on the quotient
can degenerate.  We will begin with the general case and then specialise to
torus quotients of flat space.

Let us consider an action of a Lie group \( G \) on a manifold of
dimension~\( 4d \) preserving a hypersymplectic structure and admitting a
\( G \)-equivariant moment map \( \mu \colon M \to \lie g^* \otimes \mathbb
R^3 \).  Thus, for each element~\( X \) of~\( \lie g \), the associated
component~\( \mu^X \) of~\( \mu \) satisfies~\eqref{eq:mm}.  Write \(
\mathcal G \) for the distribution on~\( M \) generated by tangent vectors
to the group action.

\begin{definition}
  \label{def:free}
  The action satisfies \emph{condition~(F)} if \( G \) acts properly and
  freely on~\( \mu^{-1}(0) \).
\end{definition}

\begin{definition}
  \label{def:smooth}
  The action satisfies \emph{condition~(S)} if at each point \( p \in
  \mu^{-1}(0) \) there is no non-zero solution to the equation
  \begin{equation}
    \label{eq:rank}
    (IX_1 + SX_2 + TX_3)_p = 0,
  \end{equation}
  with \( X_1, X_2, X_3 \in \lie g \).
\end{definition}

\begin{theorem}
  \label{thm:maxrank}
  If the \( G \)-action satisfies conditions (F) and~(S), then the
  hypersymplectic quotient \( \mu^{-1}(0)/G \) is smooth.
\end{theorem}

\begin{proof}
  It is sufficient, by (F), to show that \( \mu^{-1}(0) \) is a smooth
  manifold.  By~\eqref{eq:mm}, the kernel \( \ker d \mu \) is just the
  orthogonal complement with respect to~\( g \) of the space~\( U \)
  spanned by \( I \mathcal G \), \( S \mathcal G \) and~\( T \mathcal G \).

  As \( g \) is non-degenerate we have
  \begin{equation*}
    \rank  d \mu = 4d - \dim \ker d \mu = 4d - \dim U^{\perp} = \dim U.
  \end{equation*}
  Condition (S) implies that \( \dim U \) is \( 3 \dim \mathcal G \), which
  is \( 3 \dim G \) by~(F).  We deduce that \( d \mu \) has maximal rank
  and the result follows.
\end{proof}

\begin{remark}
  As the moment map \( \mu \) is equivariant, \( \mathcal G \) lies in \(
  \ker d \mu \) on~\( \mu^{-1}(0) \).  This implies that \( \mathcal G \)
  is orthogonal to \( I \mathcal G \), \( S \mathcal G \) and \( T \mathcal
  G \).  It follows that these spaces are mutually orthogonal.

  In the hyperk\"ahler case, where \( g \) is positive definite, this of
  course means that condition~(S) and the conclusion of the theorem follow
  automatically from the freeness of the action of \( G \), as
  in~\cite{HKLR}.
\end{remark}

It is proved in~\cite{H} that the kernels of the symplectic forms on \(
\mu^{-1}(0) \) are given by
\begin{equation}
  \label{eq:ker-I}
  \ker i^*\omega_I
  = \mathcal G
  + S (\mathcal G \cap \mathcal G^{\perp})
  + T (\mathcal G \cap \mathcal G^{\perp})
\end{equation}
and cyclically, where \( i \colon \mu^{-1}(0) \to M \) is inclusion.

\begin{definition}
  \label{def:non-deg}
  The action satisfies \emph{condition~(D)} if \( \mathcal G \cap \mathcal
  G^\bot = \{0\} \) on \( \mu^{-1}(0) \).
\end{definition}

\begin{theorem}
  \label{thm:non-degenerate}
  If the \( G \)-action satisfies conditions (F) and~(D), then the
  quotient~\( \mu^{-1}(0)/G \) inherits a smooth, non-degenerate
  hypersymplectic structure.

  On the other hand, if the \( G \)-action fulfils conditions (F) and~(S),
  then the smooth manifold~\( \mu^{-1}(0)/G \) inherits a non-degenerate
  hypersymplectic structure only if condition~(D) is satisfied.
\end{theorem}

\begin{proof}
  Hitchin's results provide a non-degenerate hypersymplectic structure
  provided \( \mu^{-1}(0)/G \) is a smooth manifold and (D)~holds: the
  symplectic form \( \omega_I' \) on the quotient is defined by the
  equation \( \pi^*\omega_I' = i^*\omega_I \), where \( \pi \) is
  projection from \( \mu^{-1}(0) \) to the quotient.  Thus for the first
  part of the Theorem we only need to consider condition~(S) of
  Definition~\ref{def:smooth}.

  Suppose \( IX_1 + SX_2 + TX_3 = 0 \) at \( p \in \mu^{-1}(0) \).  Then,
  for each \( Y \in \mathcal G_p \) we have
  \begin{equation*}
    g(X_1,Y) = g(IX_1,IY) = g(IX_1 + SX_2 + TX_3, IY) = 0
  \end{equation*}
  at~\( p \).  Thus by condition~(D), \( X_1 = 0 \).  Similarly, \( X_2 = 0
  = X_3 \) and condition~(S) holds.

  For the second part, if condition~(D) fails, then there is a \( p \in
  \mu^{-1}(0) \) such that \( V := (\mathcal G \cap \mathcal G^\bot)_p \)
  is non-zero.  Suppose \( \omega_I' \)~is non-degenerate.  Then \( \ker
  i^*\omega_I = \ker \pi_*\omega_I' = \mathcal G \).
  Equation~\eqref{eq:ker-I} gives that \( SV \) and~\( TV \) are subspaces
  of~\( \mathcal G \).  However, \( SV \) is contained in \( S\mathcal G \)
  which is orthogonal to \( \mathcal G \).  Hence, \( V \) is invariant
  under \( S \) and \( T \) and hence under \( I = ST \).  Taking \( X_1 =
  IX \), \( X_2 = SX \), \( X_3 = 0 \) for some non-zero~\( X \) in~\( V \)
  , we see that condition~(S) is violated.  Thus if the quotient geometry
  is non-degenerate and (S) holds, then (D) must hold too.
\end{proof}

As simple case of the above result is:

\begin{corollary}
  \label{cor:circ}
  A hypersymplectic quotient by a free circle action with Killing field \(
  X \) is a smooth hypersymplectic manifold except at points where \(
  g(X,X)=0 \).
  \qed
\end{corollary}

It will be useful, in the light of Corollary~\ref{cor:circ}, to have a
formula for the length of the Killing field of a circle action on \(
\mathbb C^{d,d} \).  If the action is given by
\begin{equation*}
  (z_k,w_k) \longmapsto (e^{\ii \theta_k t}z_k, e^{\ii \theta_k t}w_k),
\end{equation*}
then the associated vector field is
\begin{equation*}
  X = \sum_{k=1}^d
 \ii \theta_kz_k \frac\partial{\partial z_k} + \ii \theta_kw_k
  \frac\partial{\partial w_k} - \ii \theta_k\bar z_k\frac\partial{\partial
 \bar z_k} - \ii \theta_k\bar w_k \frac\partial{\partial \bar w_k},
\end{equation*}
giving
\begin{equation}
  \label{eq:degeneracy}
  g(X,X) = \sum_{k=1}^d \theta_k^2(\abs{z_k}^2-\abs{w_k}^2).
\end{equation}

We next investigate when the condition~(S) of Definition~\ref{def:smooth}
holds for general toric reductions of flat space.  For three vector fields
\( X_1 \), \( X_2 \), \( X_3 \) with coefficients \( \theta_k^{(i)} \), \(
i=1,2,3 \), the equation~\eqref{eq:rank} becomes
\begin{equation}
  \label{eq:rank-toric}
  \begin{aligned}
    \theta_k^{(1)}z_k &= (\ii\theta_k^{(2)}-\theta_k^{(3)})w_k, &\quad&
    \text{for \( k=1,\dots, d \),}\\
    \theta_k^{(1)}w_k &= -(\ii\theta_k^{(2)}+\theta_k^{(3)})z_k, &\quad&
    \text{for \( k=1,\dots, d \).}
  \end{aligned}
\end{equation}
For each \( k \), these equations have the form
\begin{equation*}
  az=bw,\quad aw = \bar bz,
\end{equation*}
with \( a\in \mathbb R \) and \( b\in \mathbb C \). We deduce that \( a^2z
= abw = \abs b^2z \), and similarly \( a^2 w = \abs b^2 w \). So the
system~\eqref{eq:rank-toric} has a solution if only if for each \( k \),
\begin{equation}
  \label{eq:t-rank}
  \begin{aligned}
    \text{either}\quad&(z_k,w_k)=0,\\
    \text{or}\quad&\theta_k^{(1)}=0=\theta_k^{(2)}=\theta_k^{(3)},\\
    \text{or}\quad&z_k = \xi_k w_k \quad\text{and}\quad \xi_k\theta_k^{(1)}
    = \ii (\theta_k^{(2)}+\ii\theta_k^{(3)}),\\[-0.9ex]
    &\qquad\text{for some \( \xi_k \) with \( \abs{\xi_k}=1 \).  }
  \end{aligned}
\end{equation}

We also discuss the question of degeneracy of the hypersymplectic
structure.  We know from Theorem~\ref{thm:non-degenerate} and the
discussion leading up to equation~\eqref{eq:degeneracy} that degeneracy
occurs if and only if the inner product
\begin{equation}
  \label{eq:quadratic-form}
  q = \diag ( \abs{z_k}^2 - \abs{w_k}^2 )_{k=1}^d
\end{equation}
is degenerate on \( \lie n \leqslant \mathbb R^d \) at some point \( (z,w )
\in \mu^{-1}(0) \).  Equivalently, in the notation of \S\ref{sec:moment},
at some point of \(\mu^{-1}(0) \) there is \( \zeta \in \lie n \setminus
\{0\} \) such that \( q(\zeta,\cdot ) \in \ker \iota^* = \image \beta^*\).
From~\eqref{eq:bstar}, this condition is equivalent to the existence of \(
\zeta \in \lie n \setminus \{0\} \) and \( s \in \mathbb R^n \) such that
\begin{equation}
  \label{eq:degen}
  \zeta_k ( \abs{z_k}^2 - \abs{w_k}^2) = \inp s{u_k}, \qquad \text{for
  \(k = 1, \dots, d\)}.
\end{equation}
In \S\ref{sec:quotient} we shall refine both of these criteria.

\section{Toric geometry of the quotient}
\label{sec:quotient}
We shall now study some properties of the hypersymplectic quotient~\( M \)
of~\( \mathbb C^{d,d} \) by a compact Abelian group~\( N \subset \mathbb
T^d \).  This quotient carries an action of the torus \( \mathbb T^n =
\mathbb T^d / N \) and we may consider the map \( \phi \colon M \rightarrow
\mathbb R^{3n} \) given by
\begin{equation}
  \label{eq:phidef}
  \phi \colon (z, w) \longmapsto (a,b),
\end{equation}
where \( a \) and \( b \) are as in \eqref{eq:a} and~\eqref{eq:b}.  When \(
M \)~is smooth and hypersymplectic this is the moment map for the action of
\( \mathbb T^n \) on~\( M \).

A similar map is considered in the K\"ahler and hyperk\"ahler
cases~\cite{Gu,BD}.  In the first case, the analogue of \( \phi \) is a map
\( M \rightarrow \mathbb R^n \) with image the polyhedron defined
by~\eqref{eq:poly}.  The map induces a homeomorphism of \( M/ \mathbb T^n
\) onto the polyhedron.

In the hyperk\"ahler case, one has a map onto the whole of~\( \mathbb
R^{3n} \), and again a homeomorphism \( M/\mathbb T^n \cong \mathbb R^{3n}
\).  Essentially, this follows from the \( N =\{ 1 \} \) case, i.e., the
fact that the moment map \( (z_k, w_k) \longmapsto
\left(\tfrac12(\abs{z_k}^2 - \abs{w_k}^2), \ii z_k w_k \right) \) for the
hyperk\"ahler action of \( \mathbb T^d \) on \( \mathbb H^d \) induces a
homeomorphism from \( \mathbb H^d / \mathbb T^d \) onto \( \mathbb R^{3d}
\).

In our case, the image of~\( \phi \) is an interesting subset of \( \mathbb
R^{3n} \), but we no longer obtain a homeomorphism.  Indeed the fibres
of~\( \phi \) may be disconnected.

We introduce the following notation:
\begin{equation}
  \label{eq:phi-components}
  \begin{gathered}
    a_k := \inp a{u_k} - \lambda^{(1)}_k,\\
    b_k := \inp b{u_k} - \lambda^{(c)}_k,
  \end{gathered}
\end{equation}
for \( a \in \mathbb R^n \), \( b \in \mathbb C^n \), \( k=1,\dots,d \) and
where
\begin{equation*}
  \lambda^{(c)}_k := \lambda^{(2)}_k + \ii \lambda^{(3)}_k.
\end{equation*}

\begin{proposition}
  \label{prop:phi}
  The image of the moment map~\( \phi \), equation~\eqref{eq:phidef}, is
  the set
  \begin{equation*}
    K = \bigl\{\, (a, b) \in \mathbb R^n \times \mathbb C^n : a_k
    \geqslant \abs{b_k},\ \text{for \( k=1, \dots, d\)} \,\bigr\}.
    \end{equation*}
  Moreover \( \phi \) induces a finite-to-one map \( \tilde{\phi} \) from
  \( M/\mathbb T^n \) onto \( K \).  The fibre of \( \tilde{\phi} \) over
  \( (a,b) \) has \( 2^m \) points, where \( m \) is the number of the
  inequalities in the definition of \( K \) which are \emph{strict} for \(
  (a, b) \).
\end{proposition}

\begin{proof}
  By \eqref{eq:a} and~\eqref{eq:b}, \( (z,w) \) is in \( \mu^{-1}(0)
  \subset \mathbb C^{d,d} \) and satisfies \( \phi(z,w)=(a,b) \) if and
  only if
  \begin{equation*}
    \abs{z_k}^2 + \abs{w_k}^2 = 2a_k
    \quad\text{and}\quad
    z_k\bar w_k = -\ii b_k,
  \end{equation*}
  for \( k=1,\dots,d \).  The torus \( \mathbb T^d \) acts by \(
  (z_k,w_k)\mapsto(e^{\ii\theta_k}z_k,e^{\ii\theta_k}w_k) \).  On \(
  \mu^{-1}(0) \) we have that \( \abs{w_k} \) uniquely determines \(
  (z_k,w_k) \) up to the action of~\( \mathbb T^d \).  Taking the absolute value of
  the second equation and using the first equation to eliminate \(
  \abs{z_k}^2 \), we get
  \begin{equation*}
    \abs{w_k}^2 = a_k \pm \sqrt{{a_k}^2-\abs{b_k}^2}.
  \end{equation*}
  Thus there is a solution for \( \abs{w_k}^2 \) only if \( a_k\geqslant
  \abs{b_k} \).  There are two solutions if the inequality is strict,
  otherwise there is only one solution.  The result follows.
\end{proof}

\begin{remark}
  Taking \( d=1 \) and \( N=\{1\} \) in the proof of the proposition, \(
  \phi \)~is the hypersymplectic moment map for the action of~\( \mathbb
  T^1 \) on~\( \mathbb C^{1,1} \).  One can see that the corresponding
  hypersymplectic quotient may be two points, a single point, or empty
  depending on the choice of level set.  So we may get smooth quotient sets
  of different topology for different choices of level set, contrasting
  with the hyperk\"ahler case.
\end{remark}

\begin{proposition}
  The set \( K = \phi(M) \) is convex in \( \mathbb R^{3n} \).
\end{proposition}

\begin{proof}
  Using~\eqref{eq:phi-components}, the set \( K \) is the intersection of
  sets
  \begin{equation*}
    K_k = \{\, (a,b) \in \mathbb R^n \times \mathbb C^n :
    a_k\geqslant\abs{b_k} \,\}, \qquad k=1,\dots,d.
  \end{equation*}
  However, \( K_k \)~is the preimage of the solid cone \( \{\,(x,z)\in
  \mathbb R\times \mathbb C: x\geqslant\abs z\,\} \) under the affine map
  \( (a,b) \mapsto (a_k,b_k) \).  Thus \( K_k \) is convex, and it follows
  that \( K \) is convex too.
\end{proof}

\begin{corollary}
  \( M \)~is connected if and only if each inequality in the definition
  of~\( K \) is an equality at some point of~\( K \).  \qed
\end{corollary}

\begin{theorem}
  \label{thm:compact}
  The hypersymplectic quotient \( M = \mu^{-1}(0)/N \) is compact if and
  only if the vectors \(u_1, \dots, u_d \) define a \emph{bounded}
  polyhedron in \( \mathbb R^n \).
\end{theorem}

\begin{proof}
  Suppose \( u_1,\dots, u_d \) define a bounded polyhedron.  As \( M =
  (\mu_I^{-1}(0) \cap (\mu_S + \ii \mu_T)^{-1} (0))/ N \), it is enough to
  show compactness of \( \mu_I^{-1}(0) \).

  If \( (z,w) \in \mu_I^{-1}(0) \), then the vector \( a \) of \eqref{eq:a}
  must live in the polyhedron \( \{\, a : \inp a{u_k} \geqslant
  \lambda_k^{(1)} \,\} \), which is compact by hypothesis.  Now
  \eqref{eq:a} gives us a bound on the \( \abs{z_k} \) and~\( \abs{w_k} \)
  in terms of the \(u_k \) and~\(\lambda_k \).

  Conversely, note that \( M \) is compact if and only if \( M/\mathbb T^n
  \) is compact, and hence if and only if \( K \) is compact.  We may
  define a projection \( p \colon K \rightarrow \mathbb C^n \) by \( p
  (a,b) = b \).  The fibres of~\( p \) are the polyhedra
  \begin{equation*}
    F_b = \bigl\{\, a \in \mathbb R^n :  \inp a{u_k} \geqslant
      \lambda_k^{(1)} + \abs{b_k} \,\bigr\}.
  \end{equation*}
  If \( M \) is compact then \( F_b \) is compact for each \( b \in p(K)
  \).  As \( F_b \) is non-empty, this implies that \(u_1, \dots, u_d \)
  define a bounded polyhedron as in~\eqref{eq:poly}.
\end{proof}

As promised in~\S\ref{sec:non-degeneracy}, we now refine the criteria for
smoothness and non-degeneracy of the hypersymplectic quotient~\( M \),
in terms of the map~\( \phi \) and the vectors~\( u_k \).

\begin{figure}[tbp]
  \centering
  \includegraphics{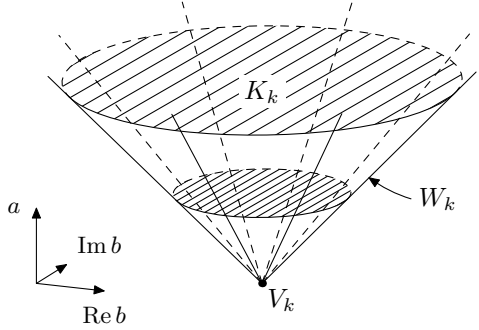}
  \caption{The solid cone \( K_k \) when \( n=1 \).}
  \label{fig:Kk}
\end{figure}

In order to consider smoothness of~\( M \), we need to discuss the orbit
types of the actions of~\( N \) on~\( \mu^{-1}(0) \) and the smoothness of
\( \mu^{-1}(0) \) itself.  Let us begin with the orbit types.  In the
notation~\eqref{eq:phi-components}, put
\begin{equation}
  \label{eq:vertex}
  V_k :=  \{\, (a,b) \in \mathbb R^n\times\mathbb C^n : a_k=0=b_k
  \,\},
\end{equation}
see Figure~\ref{fig:Kk}.  Note that \eqref{eq:a} and~\eqref{eq:b} show that
for \( (z,w) \in \mu^{-1}(0) \)
\begin{equation*}
  z_k=w_k=0 \qquad  \text{if and only if} \qquad (a,b)=\phi(z,w) \in V_k.
\end{equation*}
Following~\cite{Gu}, if \( A \) is a subset of \( \{1, \dots, d \} \) we
denote by \( \mathbb T_A \) the torus whose Lie algebra is
\begin{equation*}
  \mathbb R_A = \Span_{\mathbb R} \{\, e_k : k \in A \,\}.
\end{equation*}
We deduce that the stabiliser of \( (z,w) \), for the \( \mathbb T^d \)
action, is \( \mathbb T_J \), where
\begin{equation}
  \label{eq:J}
  J := \{\, k : \phi(z,w) \in V_k \,\}.
\end{equation}

\begin{proposition}
  \label{prop:stabiliser}
  For \( (z,w) \in \mu^{-1}(0) \):
  \begin{enumerate}
  \item[(i)] \( \Stab_N (z,w) \) equals \( \mathbb T_J \cap N \), where \(
    \mathbb T_J \) is the torus whose Lie algebra is spanned by the vectors
    \( e_k \) for which \( \phi(z,w)\in V_k \);
  \item[(ii)] On \( M \), \( \Stab_{\mathbb T^n} (z,w) \) is the torus
    whose Lie algebra is spanned by the vectors \( u_k \) such that \(
    \phi(z,w)\in V_k \).
  \end{enumerate}
  Moreover, putting \( \mathcal B_x = (u_k : x \in V_k) \),
  \begin{enumerate}
  \item[(iii)] \( \Stab_N(z,w) \) is finite for all \( (z,w) \in
    \mu^{-1}(0) \) if and only if for each \( x \in \phi(M) \), the
    collection of vectors~\( \mathcal B_x \) is linearly independent,
  \item[(iv)] \( \Stab_N(z,w) = 1 \) for all \( (z,w) \in \mu^{-1}(0) \) if
    and only if for each \( x \in \phi(M) \), the collection of vectors~\(
    \mathcal B_x \) is contained in a \( \mathbb Z \)-basis for \( \mathbb
    Z^n \).
  \end{enumerate}
\end{proposition}

\begin{proof}
  Statements (i) and (ii) follow from the above discussion, while (iii) and
  (iv) follow from (i), (ii) and results of Delzant~\cite{De} and
  Guillemin~\cite{Gu}, as cited in the proof of Theorems 3.2 and 3.3
  in~\cite{BD}.  To see this, note that \( V_k \) is the affine flat \( H_k
  \) introduced in~\cite{BD}.
\end{proof}

\begin{remark}
  The argument of Theorems 3.2 and 3.3 of~\cite{BD} shows that if every \(
  n+1 \) of the \( V_k \) have empty intersection, then the condition
  of~(iii) holds and hence \( N \) acts locally freely on \(\mu^{-1}(0) \).
  In particular, for any given collection of vectors \( u_k \), the action
  will be locally free for generic choice of \( \lambda \).  Hence, by
  Sard's theorem, \(M = \mu^{-1}(0)/N \) will have at worst orbifold
  singularities for generic choice of \( \lambda \).

  Similarly, we see that if for each \( A \subset \{1,\dots,d\} \) of
  size~\( n \), the collection \( (u_k : k\in A) \) is a \( \mathbb Z
  \)-basis for \( \mathbb Z^n \), then \( M = \mu^{-1}(0)/N \) is a
  manifold for generic~\( \lambda \).
\end{remark}

To determine precise conditions for the smoothness of~\( \mu^{-1}(0) \),
let
\begin{equation*}
  W_k := \{\, (a,b)\in \mathbb R^n \times \mathbb C^n : a_k =
  \abs{b_k} \,\}.
\end{equation*}
A point \( (z,w) \in \mu^{-1}(0) \) has \( \phi(z,w) \in W_k \) if and only
if \( \abs{z_k}=\abs{w_k} \).  Note that \( V_k\subset W_k \) and that both
sets may be empty for a given~\( k \).  Let \( J \) be as in~\eqref{eq:J}
above and put
\begin{equation*}
  L := \{\, \ell : \phi(z,w)\in W_\ell \,\}.
\end{equation*}
Note that \( L \) contains~\( J \).

Consider the equations~\eqref{eq:t-rank} at \( (z,w) \in \mu^{-1}(0) \).
For \( k\in J \), there is no restriction on \( \theta^{(i)}_k \).  For \(
k\in L' \), the complement of \( L \) in \( \{1,\dots,d\} \), we require \(
\theta^{(i)}_k=0 \), for \( i=1,2,3 \).  For \( k\in L\setminus J \), we
have \( \xi_k\theta^{(1)}_k = \ii(\theta^{(2)}_k+\ii\theta^{(3)}_k)) \).
For such~\( k \), we have \( \abs{z_k}=\abs{w_k}\ne0 \), so \( \xi_k \)~is
uniquely determined by
\begin{equation*}
  \xi_k = \frac{z_k}{w_k} = \frac{z_k\bar w_k}{\abs{w_k}^2} = -\ii
  \frac{b_k}{a_k}.
\end{equation*}
Consider the map
\begin{equation}
  \label{eq:Lambda}
  \begin{gathered}
    \Lambda_{(a,b)} \colon \lie n_{L,J}\otimes (\mathbb R\times
    \mathbb C) \to \mathbb C_{L\setminus J},\\
    \Lambda_{(a,b)} (\theta^{(1)}_k,\theta^{(c)}_k) = (b_k \theta^{(1)}_k +
    a_k \theta^{(c)}_k),
  \end{gathered}
\end{equation}
where \( \lie n_{L,J} \) is the projection to \( \mathbb R_{J'} \) of \(
\lie n_L := \lie n \cap \mathbb R_L \).  Alternatively, \( \lie n_{L,J} \)
is the kernel of the map \( \beta_{L,J}\colon \mathbb R_{L\setminus J} \to
\mathbb R^n / \image(\beta|_{\mathbb R_J}) \) induced by~\( \beta \).

\begin{proposition}
  \label{prop:smoothness}
  For a locally free action of~\( N \) on~\( \mu^{-1}(0) \), the
  condition~(S) of Definition~\ref{def:smooth} for smoothness of the level
  set~\( \mu^{-1}(0) \) holds if and only if the linear map~\(
  \Lambda_{(a,b)} \) of equation~\eqref{eq:Lambda} is injective at each
  point \( (a,b) \) of~\( K \).  \qed
\end{proposition}

\begin{remark}
  Over the \emph{combinatorial interior} of~\( K \),
  \begin{equation*}
    \CInt(K) = K \setminus \bigcup_{k=1}^d W_k,
  \end{equation*}
  the set~\( L \) is empty, so \( \phi^{-1}(\CInt K) \subset M \) is always
  a smooth manifold.  The combinatorial interior is equal to the
  topological interior if \( K_k \ne W_k \), for all~\( k \).
\end{remark}

\begin{remark}
  \label{rem:smooth-bound}
  By Proposition~\ref{prop:stabiliser}, if the action of~\( N \) is locally
  free then \( \lie n_L \) is transverse to~\( \mathbb R_J \) and hence the
  orthogonal projection \( \lie n_L \to \lie n_{L,J} \) is injective.  Thus
  the smoothness condition will necessarily fail if \( 3\dim \lie n_L >
  2\abs{L\setminus J} \).  Since \( \dim \lie n_L \geqslant \abs{L} - n \),
  we conclude that smoothness requires \(\abs L \leqslant 3n \), i.e., no
  more than \( 3n \) of the \( W_k \)'s may meet in~\( K \).
\end{remark}

Finally, let us consider non-degeneracy of the hypersymplectic structure.

\begin{proposition}
  \label{prop:phi-non-deg}
  The non-degeneracy condition~(D) of Definition~\ref{def:non-deg} fails at
  some point \( (z,w) \in \mu^{-1}(0) \) with \( \phi(z,w) = (a,b) \) if
  and only if there exist scalars \( \zeta_1, \ldots, \zeta_d \), not all
  zero, and \( s \in \mathbb R^n \) such that
  \begin{equation}
    \label{eq:proportionality}
    4{\zeta_k}^2 \left( {a_k}^2 - \abs{b_k}^2\right) = {\inp s{u_k}}^2,
    \qquad \text{for \(k=1, \ldots, d \)}
  \end{equation}
  and
  \begin{equation}
    \label{eq:lin-dep}
    \sum_{k=1}^d \zeta_k u_k =0.
  \end{equation}
\end{proposition}

\begin{proof}
  This is immediate from the the discussion at the end
  of~\S\ref{sec:non-degeneracy}, equation~\eqref{eq:degen}, the definition
  of \( \lie n \) as~\( \ker \beta \), and the proof of
  Proposition~\ref{prop:phi} which expresses \( \abs{z_k} \) and \(
  \abs{w_k} \) in terms of~\( a \) and \( b \).
\end{proof}

\begin{remark}
  \label{rem:s-zero}
  Consider the special case when \( s = 0 \).  The expression \( {a_k}^2 -
  \abs{b_k}^2 \) is zero on~\( W_k \). Suppose \( A \subset \{1,\dots,d\}
  \) has \( n+1 \) elements.  Then \( (u_k : k \in A) \) is linearly
  dependent, so we may find \( (\zeta_1, \dots , \zeta_d ) \in \mathbb R_A
  \setminus \{0\} \) satisfying~\eqref{eq:lin-dep}.  If
  \begin{equation*}
    K \cap \bigcap_{k\in A} W_k \ne \varnothing
  \end{equation*}
  then we then obtain a solution to the remaining
  equations~\eqref{eq:proportionality}. Thus if \( n+1 \) of the \( W_k
  \)'s meet in~\( K \), then condition~(D) necessarily fails.  This is a
  considerably stronger restriction than that obtained for smoothness in
  Remark~\ref{rem:smooth-bound}, and as we will see
  in~\S\ref{sec:examples}, one may easily obtain smooth quotients with
  hypersymplectic structures that degenerate on some hypersurface.
\end{remark}

\begin{theorem}
  Let \( \phi \) be the moment map for the action of \( \mathbb T^n \) on
  the hypersymplectic quotient \( M = \mu^{-1}(0)/N \).  Suppose the
  combinatorial interior \( \CInt(K) \) of the image~\( K \) of~\( \phi \)
  is non-empty.  Then an open subset of \( M \) carries a smooth
  non-degenerate hypersymplectic structure.  
  
  When \( M \) is compact, the degeneracy locus is non-empty of codimension
  at least one.
\end{theorem}

Non-trivial examples of quotients with empty degeneracy locus will be given
in~\S\ref{sec:examples}.

\begin{proof}
  We have already noted that \( M \) is smooth over~\( \CInt(K) \).  For \(
  (a,b) \in \CInt(K) \), we have \( a_k^2\ne\abs{b_k}^2 \) for all~\( k
  \).  Equations~\eqref{eq:proportionality} may be thus be solved for \(
  \zeta_k \) and equation~\eqref{eq:lin-dep} becomes
  \begin{equation*}
    H(s,w) := \sum_{k=1}^d w_k \inp s{u_k}u_k = 0,
  \end{equation*}
  with \( w=(w_1,\dots,w_d) \), \( w_k =
  1/(2\varepsilon_k\sqrt{{a_k}^2-\abs{b_k}^2}) \), \( \varepsilon_k=\pm 1
  \) and \( s\ne0 \).  Now \( H\colon R_1\times R_2 \to \mathbb R^n \),
  where \( R_1 = \mathbb R^n\setminus\{0\} \) and \( R_2 \) is an \( n
  \)-dimensional manifold contained in the complement of the coordinate
  axes in~\( \mathbb R^d \), has differential \(
  DH_{(s,w)}(\varsigma,\upsilon) = \sum_{k=1}^d
  \inp{w_k\varsigma+\upsilon_k s}{u_k}u_k \) which has rank \( n \) when \(
  w_k\ne0 \) for all~\( k \).  Thus \( H^{-1}(0) \) is a submanifold of \(
  R_1\times R_2 \subset \mathbb R^n\times\mathbb R^d \) of dimension~\( n
  \).  But if \( (s,w) \) lies in \( H^{-1}(0) \) then so does \( (\lambda
  s,w) \) for \( \lambda\in\mathbb R\setminus{0} \).  So the projection of
  \( H^{-1}(0) \) to the second factor \( R_2 \), has dimension at most \(
  n-1 \), i.e., the degeneracy locus is at least codimension one.
  
  Now suppose that \( M \) is compact.  Consider the map \( p\colon
  K\to\mathbb C^n \) of Theorem~\ref{thm:compact} and its fibres \( F_b \).
  Note that when the interior \( \Int F_b \) is non-empty we have \( \Int
  F_b = \{\,a\in\mathbb R^n: a_k>\abs{b_k}\,\} \).  As this is an open
  condition on~\( b \), we see that \( p(\Int F_b) \subset \Int(p(K)) \).
  Fix \( b \) on the boundary of \( p(K) \).  Then \( F_b \) is a compact
  polytope in \( \mathbb R^n \) with empty interior.  Let \( v \) be a
  vertex of that polytope.  Then \( v \) is at the intersection of at least
  \( n \) hyperplanes \( a_k=\abs{b_k} \).  However, if only \( n \)
  hyperplanes meet, then we can find interior points of \( F_b \) close to
  \( v \).  Thus \( n+1 \) hyperplanes meet in~\( v \).  This is the same
  as saying that \( (v,b) \) lies on \( n+1 \) of the~\( W_k \).  By
  Remark~\ref{rem:s-zero}, this implies that the hypersymplectic structure
  is degenerate at \( (v,b) \).
\end{proof}

We shall next relate our quotients to toric varieties.  We first prove a
lemma.

\begin{lemma}
  \label{lem:compact}
  Let \( \delta \colon \mathbb T^d \to \mathbb T^d \times \mathbb T^d =
  \mathbb T^{2d} \) be the diagonal map, and let \( N \) be a compact
  Abelian subgroup of~\( \mathbb T^d \).

  If \( N \leqslant \mathbb T^d \) is defined by a collection of vectors in
  \(\mathbb R^n \) defining a bounded polyhedron, then \( \delta(N)
  \leqslant \mathbb T^{2d} \) is defined by a collection of vectors in \(
  \mathbb R^{n+d} \) defining a bounded polyhedron.
\end{lemma}

\begin{proof}
  We have an exact sequence
  \begin{equation*}
    0 \longrightarrow \delta_* (\lie n )
    \overset{\tilde{\iota}}{\longrightarrow}
    \mathbb R^d \oplus \mathbb R^d = \mathbb R^{2d}
    \overset{\tilde{\beta}}{\longrightarrow}
    \mathbb R^{n+d} \longrightarrow 0,
  \end{equation*}
  for some map \( \tilde{\beta} \).  It is straightforward to check that
  \(\delta_*(\lie n) = \{\,(v,v) : v \in \lie n \,\} \) is the kernel of the
  map defined by
  \begin{equation*}
    e_k \longmapsto
    \tilde u_k :=
    \begin{cases}
      u_k + e_{k+d},&\text{if \( k\leqslant d \),}\\
      -e_k,&\text{if \( k > d \).}
    \end{cases}
  \end{equation*}
  So we can take this map to be \( \tilde{\beta} \), and hence
  \(\delta_*(\lie n) \) is defined by vectors \( \tilde u_1,\dots,\tilde
  u_{2d} \).

  Consider the polyhedron \( \{\, s \in \mathbb R^{n+d} : \inp s{\tilde
  u_k} \geqslant c_k,\ k=1, \dots, 2d \,\} \).  The defining inequalities
  may be written as
  \begin{alignat*}{2}
    \inp{s_1}{u_k} + \inp{s_2}{e_{k+d}} &\geqslant c_k,&\qquad&
    \text{for \( k=1, \dots, d \),} \\
    \inp{s_2}{-e_{k+d}} &\geqslant c_{k+d},&&\text{for \( k=1, \dots, d \),}
  \end{alignat*}
  where \( s=s_1 + s_2 \), \( s_1 \in \mathbb R^n \) and \(
  s_2 \in \Span \{e_{d+1}, \dots, e_{2d} \} \).  We see that \(
  \inp{s_1}{u_k} \geqslant c_k + c_{k+d}\), for \(k=1, \dots, d \).  As \(
  u_1, \dots, u_d \) define a bounded polyhedron by hypothesis, we get a
  bound on \( s_1 \), and hence, from
  \begin{equation*}
    c_k - \inp{s_1}{u_k}  \leqslant \inp {s_2}{e_{k+d}} \leqslant
    -c_{k+d}
  \end{equation*}
  a bound on \(s_2 \).
\end{proof}

Now, observe that \(\mu_I \) is just the moment map \( m \) for the
K\"ahler action of \(\delta(N)\) on \(\mathbb C^{2d} \).  (Here \( \mathbb
C^{2d} \) has the standard complex structure \( I_0 \), not the
hypersymplectic complex structure \( I \), as remarked
in~\S\ref{sec:flat}).

The construction of Guillemin and Delzant now gives:

\begin{theorem}
\label{thm:toricv}
  \( M \) is the sub-variety
  \begin{equation*}
    \mu_S + \ii \mu_T =0
  \end{equation*}
  in the toric variety \( m^{-1}(0)/\delta(N) \).  If the vectors \( u_k \)
  define a bounded polyhedron, then by Lemma~\ref{lem:compact} this toric
  variety is compact. \qed
\end{theorem}

\begin{remark}
  Note that \( \mu_S + \ii \mu_T \), although \( I \)-holomorphic, is
  \emph{not} holomorphic with respect to the complex structure on this
  toric variety, which is induced from \( I_0 \).
\end{remark}

We conclude by discussing two actions that can occur on \( M \) for special
choices of \( \lambda \).

\begin{remark}
  If we take \( \lambda_k^{(j)}=0 \) for all \( j,k \), then \( (0,0) \in
  \mu^{-1}(0) \) and is a fixed point of \( N \), giving a singular point
  in the quotient \( M = \mu^{-1}(0)/N \).  In fact it follows from
  \eqref{eq:a}, \eqref{eq:b} with \( \lambda_j^{(k)}=0 \) that we have a
  scaling action \( (z, w) \mapsto (t.z, t.w) \), for \( t\in \mathbb R^*
  \) on \( \mu^{-1}(0) \) which descends to the quotient, so that \( M \)
  is a cone with vertex at \( (0,0) \).
\end{remark}

\begin{remark}
  \label{rem:circleaction}
  Harada and Proudfoot~\cite{HP} have observed that in the hyperk\"ahler
  case, if we take \( \lambda_k^{(2)} + \ii \lambda_k^{(3)}=0 \) for all \(
  k \), then the quotient \( M \) admits a circle action
  \begin{equation*}
    (z, w) \mapsto (z, e^{\ii \psi} w).
  \end{equation*}
  This action is holomorphic with respect to \( I \) but not with respect
  to \( J \) or \( K \).

  This action does not occur for the K\"ahler quotients of~\cite{Gu,De}.
  However we observe that it does exist for our hypersymplectic quotients,
  provided we take \( \lambda_k^{(2)} + \ii \lambda_k^{(3)} =0 \) as above.
  The action is compatible with \( I \) but not with \( S \) or \( T \).

  Note also that under this condition on \( \lambda \), the hypersymplectic
  quotient contains two distinguished subvarieties, defined by the
  vanishing of \( z \) and \( w \) respectively.  Each of these may be
  identified with a K\"ahler quotient of \( \mathbb C^d \) by \( N \), and
  hence with a toric variety.

  The locus \(w=0 \) lies in the fixed point set of the circle
action. As \( M \) is a quotient, we may in addition have other
components of the fixed point set. We find, as in the hyperk\"ahler
case, that the fixed point set in general is a union of toric
varieties which may be enumerated in terms of conditions on the vectors
\(u_k \).

\end{remark}

\section{Involutions}
\label{sec:involutions}
In the K\"ahler case studied by Guillemin and Delzant, the moment map
\( \mu \colon \mathbb C^d \mapsto \lie n^* \) is invariant under the
involution of \( \mathbb C^d \) given by complex conjugation.  In
fact, conjugation induces an involution \( \gamma \) of the quotient
\( M = \mu^{-1}(0)/N \).  Moreover the fixed point set of complex
conjugation in \( \mu^{-1}(0) \) is a cover of the fixed point set of
\( \gamma \) in \( M \).  The group of deck transformations is the
finite group \( \Gamma \) of involutions in \( N \).

In the hyperk\"ahler situation there appears to be no such involution in
general.

In our hypersymplectic case, however, we \emph{do} have an involution.
Explicitly, the map \( \sigma \colon \mathbb C^{d,d} \mapsto
\mathbb C^{d,d} \) given by
\begin{equation*}
  \sigma \colon (z_k, w_k)  \mapsto (\bar w_k,\bar z_k)
\end{equation*}
preserves \( \mu^{-1}(0) \), and sends \( N \)-orbits to \( N \)-orbits.
(It does not commute with the action, but we have \( \sigma (g \cdot (z,w))
= g^{-1} \cdot \sigma(z,w) \), which suffices).

It follows that \( \sigma \) induces an involution \( \hat{\sigma} \) on
the hypersymplectic quotient \( M = \mu^{-1}(0)/N \).  Let us denote by \(
\mu^{-1}(0)_\sigma \) and \( M_{\hat{\sigma}} \) the fixed point sets of
the involutions in \( \mu^{-1}(0) \) and \( M \) respectively.

\begin{theorem}
  The natural projection \( \pi \colon \mu^{-1}(0) \rightarrow M \) induces
  a surjection
  \begin{equation*}
    \mu^{-1}(0)_{\sigma} \rightarrow M_{\hat{\sigma}}.
  \end{equation*}
  If \( N \) acts freely on \( \mu^{-1}(0) \) this map is a cover whose
  group of deck transformations is the finite group
  \begin{equation*}
    \Gamma = \{\, h \in N : h^2 = 1 \,\}.
  \end{equation*}
\end{theorem}

\begin{proof}
  Observe first that \( (z, w) \in \mu^{-1}(0) \) represents a point in \(
  M_{\hat{\sigma}} \) if and only if there exists \( g \in N \) with \( w =
  g^{-1}\bar z = \overline{gz} \).  Now if \( h \in N \) satisfies \( h^2 =
  g \), then
  \begin{equation*}
    h\cdot(z,w) = h\cdot(z, \overline{gz}) = (hz, \overline{h^{-1}gz})
    = (hz, \overline{hz})
  \end{equation*}
  so \( (z,w) \) represents the same point in~\( M \) as does \( (hz,
  \overline{hz}) \in \mu^{-1}(0)_{\sigma} \).  This proves the surjectivity
  assertion.

  Next, suppose that two points \( (z,\bar z) \) and \( (u, \bar u) \) in
  \( \mu^{-1}(0)_{\sigma} \) are related by the action of \( g \in N \).
  We need, for each \( k \),
  \begin{equation*}
    z_k = e^{\ii \theta_k} u_k \quad\text{and}\quad
    \bar z_k =e^{\ii \theta_k} \bar u_k
  \end{equation*}
  for some \( e^{\ii \theta_k} \in \mathbb T^1 \).  Hence either \( e^{i
  \theta_k} \) is an order two element in \( \mathbb T^1 \) or \( z_k =
  u_k=0 \).  We deduce that \( g^2 \in \mathbb T_J \cap N \) where \( J \)
  is the set of indices for which \( z_k=0 \).

  If \( N \) acts freely on \( \mu^{-1}(0) \) then \( \mathbb T_J \cap N \)
  is trivial by Proposition~\ref{prop:stabiliser}.  The remaining
  assertions now follow easily.
\end{proof}

\begin{proposition}
  The symplectic forms \( \omega_I, \omega_S \) and \( \omega_T \) all
  vanish when restricted to \( M_{\hat{\sigma}} \).
\end{proposition}

\begin{proof}
  The involution \( \hat\sigma \) pulls back \( \omega_I, \omega_S \) and
  \( \omega_T \) to their negatives.
\end{proof}

\begin{remark}
  Non-degenerate hypersymplectic manifolds with multi-Lagrangian subsets
  have recently appeared in~\cite{FPPW}.  There hypersymplectic structures
  arise on \( \mathbb T^{2n} \)-fibrations over a \( \mathbb T^{2n} \)
  base, and the fibres are multi-Lagrangian.
\end{remark}

In contrast, we have:

\begin{proposition}
  \label{prop:fixed-deg}
  The fixed point set \( M_{\hat\sigma} \) is non-empty only if
  \begin{equation*}
    \bigcap_{k=1}^d W_k \ne \varnothing.
  \end{equation*}
  The non-degeneracy statement~(D), Definition~\ref{def:non-deg}, is not
  satisfied on \( TM |_{M_{\hat\sigma}} \). In particular, if the action
  satisfies conditions~(F) and~(S) of Definitions~\ref{def:free}
  and~\ref{def:smooth}, then \( M \) is smooth but the hypersymplectic
  structure is degenerate over \( M_{\hat\sigma} \).
\end{proposition}

\begin{proof}
  If \( (z,w) \in \mu^{-1}(0) \) is a fixed point for~\( \sigma \), then \(
  \abs{z_k} = \abs{w_k} \) for \( k = 1,\dots,d \).  This is the condition
  that \( \phi(z,w) \in W_k \) for \( k=1,\dots,d \) and the assertion
  follows.

  When \( \abs{z_k} = \abs{w_k} \) for all~\( k \), the quadratic
  form~\eqref{eq:quadratic-form} is identically zero.  Thus condition~(D)
  fails.  The final statement follows from Theorems \ref{thm:maxrank}
  and~\ref{thm:non-degenerate}.
\end{proof}

\section{Examples}
\label{sec:examples}

\subsection{The diagonal circle action}
\label{sec:diagonal}

Let us take \( d = n+1 \), and take
\begin{gather*}
  u_k = e_k, \qquad\text{for \( k=1, \dots, n \),} \\
  u_{n+1} = -(e_1 + \dots + e_n)
\end{gather*}
so the vectors \( u_k \) define the standard simplex in \( \mathbb R^n \).
Now \( N \) is the standard diagonal circle in \( \mathbb T^{n+1} \).  The
corresponding toric hyperk\"ahler manifold is the Calabi space~\(
T^*\mathbb CP(n) \).

In the hypersymplectic case \( M=\mu^{-1}(0)/N \) is the quotient of the
subset of \( \mathbb C^{n+1,n+1} \) cut out by the equations
\begin{gather*}
  \sum_{k=1}^{n+1} \abs{z_k}^2 +
  \abs{w_k}^2 = -2 \sum_{k=1}^{n+1} \lambda_k^{(1)}, \\
  \sum_{k=1}^{n+1} z_k \bar w_k = -\ii\sum_{k=1}^{n+1} ( \lambda_k^{(2)}
  + \ii \lambda_k^{(3)}),
\end{gather*}
by the action
\begin{equation*}
  z_k \longmapsto e^{\ii \theta} z_k,\quad w_k \longmapsto e^{\ii
  \theta} w_k.
\end{equation*}
Using Theorem \ref{thm:toricv}, this may be identified with a hypersurface
in \( \mu_I^{-1}(0)/N = \mathbb C \mathbb P^{2n+1} \).

Note that the action
\begin{equation*}
  (z,w) \longmapsto (Az, Aw),\qquad\text{for \( A \in U(n+1) \)}
\end{equation*}
preserves the level set of \(\mu^{-1}(0) \) and commutes with the action of
\( N \), so defines an effective \( PU(n+1) \) action on \( M \).

The stabiliser of \( (z,w) \) is \(P (U(1) \times U(n)) \) if \(z,w\) are
linearly dependent, and \(P(U(1) \times U(n-1)) \) otherwise: the \(P
U(n+1) \) action on \(M \) is therefore cohomogeneity one.

Let us write
\begin{equation*}
  P = -2 \sum_{k=1}^{n+1} \lambda_k^{(1)},\qquad
  Q = -\ii\sum_{k=1}^{n+1} ( \lambda_k^{(2)} + \ii \lambda_k^{(3)})
\end{equation*}
A necessary condition for \( M \) to be nonempty is \( \abs Q \leqslant
\frac12P \)

The vector field \( X \) for the action of \( N \) takes the value \( (\ii
z, \ii w) \) at the point \( (z,w) \).  From~\S\ref{sec:non-degeneracy} we
see that \( IX, SX, TX \) are linearly independent unless \(w = \lambda z
\) for a complex number \( \lambda\) of unit modulus.  It is easy to check
that such a point \( (z,w) \) cannot lie in \( \mu^{-1}(0) \) except in the
special case \( \abs Q = \frac12 P \).  We deduce that \( M \) is a smooth
manifold unless \( \abs Q =\frac12P \).

The hypersymplectic structure on \( M \) will degenerate at some points,
however.  Observe that the point \( (z,w) \), where
\begin{equation*}
  z = (\sqrt{\tfrac12P},0, \dots,0),\qquad
  w = (\frac{\bar Q}{\sqrt{\frac12P}},w_2, \dots, w_{n+1})
\end{equation*}
and
\begin{equation*}
  \sum_{k=2}^{n+1} \abs{w_k}^2 = \tfrac12 P - \frac{\abs Q^2}{\frac12 P}
\end{equation*}
lies in \( \mu^{-1}(0) \) when this set is nonempty.  Moreover, from
\eqref{eq:degeneracy}, we have that the Killing field for the circle \( N
\) is null at \( (z,w) \).  Hence the hypersymplectic structure degenerates
at this point.

If we take \( Q =0 \) then we have a circle action
\begin{equation*}
  (z, w) \mapsto (z, e^{\ii \psi}w )
\end{equation*}
as discussed in Remark~\ref{rem:circleaction}.  The fixed point set of this
action in \( M = \mu^{-1}(0)/N \) is the union of the loci \(z=0 \) and \(
w=0 \).  These are both diffeomorphic to \( \mathbb C \mathbb P^n \), and
are the special orbits of the \( PU(n+1) \) action in this case.

\subsection{Codimension-one subgroups}
\label{sec:codim1}

Take \( n=1 \), and first consider
\begin{equation*}
  u_k = e_1,\qquad\text{for \( k=1, \dots,  d \),}
\end{equation*}
so
\begin{equation*}
  \lie n = \Bigl\{\sum_{k=1}^d a_k e_k :  \sum_{k=1}^d a_k = 0
  \Bigr\}
\end{equation*}
and \( N \) is the torus \( \mathbb T^{d-1} = \{ (t_1, \dots, t_d) : t_1t_2
\dots t_d = 1 \} \) in \( \mathbb T^d \).  The Gibbons-Hawking
multi-instanton metrics are obtained as hyperk\"ahler quotients of \(
\mathbb H^d \) by~\( N \).

We obtain, analogously, 4-dimensional hypersymplectic quotients of \(
\mathbb C^{d,d} \) by \( N \), with an action of the circle group \(
\mathbb T^d/N \).  Like the multi-instanton spaces, these hypersymplectic
spaces are non-compact, as can by seen from Theorem~\ref{thm:compact}.

Note that the sets \( V_k \) of~\eqref{eq:vertex} are now just the \( d \)
points
\begin{equation*}
  V_k = \{(\lambda_k^{(1)}, \lambda_k^{(2)},
  \lambda_k^{(3)})\},\qquad\text{for \( k=1, \dots, d \)}.
\end{equation*}
Proposition~\ref{prop:stabiliser} (iv) shows that \( N \) acts freely on
\( \mu^{-1}(0) \) provided that these \( d \) points are distinct (as in
the hyperk\"ahler case).

The set \( K = \phi(M) \) in Proposition~\ref{prop:phi} is the intersection
of the \( d \) cones
\begin{equation*}
  \{ (a, b) \in \mathbb R \times \mathbb C  :  a - \lambda_k^{(1)}
  \geqslant | b - (\lambda_k^{(2)} + \ii \lambda_k^{(3)})| \}
\end{equation*}
with vertices at \( (\lambda_k^{(1)}, \lambda_k^{(2)}, \lambda_k^{(3)}) \),
for \( k=1, \dots, d \).  All these cones have the same angle and parallel
axes.

The fixed points of the circle action correspond, from
Proposition~\ref{prop:stabiliser}(ii), to the preimages under \( \phi \) of
those points \( (\lambda_k^{(1)}, \lambda_k^{(2)}, \lambda_k^{(3)}) \)
which lie in the intersection of the cones.

We deduce that different configurations of the cones may give a zero or
non-zero number of fixed points.  This gives an example of different choices
of level set giving quotient sets which are inequivalent as \( \mathbb T^n
\)-manifolds.

By contrast, in the hyperk\"ahler case, \( \phi(M) \)~is the whole of~\(
\mathbb R^3 \), and we always get \( d \)~fixed points for the circle
action (provided that the points~\( H_k=V_k \) are distinct).

\smallskip
Let us now consider what may happen for arbitrary choices of \( u_k \) when
\( n=1 \) and \( d \)~is small.

\begin{example}
  Consider \( n=d=1 \).  This case is relatively trivial as the group~\( N
  \) is discrete.  However, it illustrates a number of features of our
  constructions.
  
  The map \( \beta \) is defined by a single \( u_1 \in \mathbb
  Z\setminus\{0\} \), which we may take to be positive.  The image of the
  moment map~\( \phi \) is the solid cone \( K=K_1 \) consisting of \(
  (a,b) \in \mathbb R\times\mathbb C \) such that \( a-\lambda_1^{(1)}/u_1
  \geqslant \abs{b-\lambda_1^{\smash{(c)}}/u_1} \), as in
  Figure~\ref{fig:Kk}.  The vertex \( V_k \) is the point \(
  (\lambda_1^{(1)}/u_1,\lambda_1^{(c)}/u_1) \).  The \( N \)-action is free
  only if \( u_1=1 \), in which case \( N=\{1\} \).  As \( \lie n=\{0\} \),
  the maps \( \Lambda_{(a,b)} \)~\eqref{eq:Lambda} are injective, so the
  quotient is indeed smooth when \( u_1=1 \).  The
  equation~\eqref{eq:lin-dep} implies \( \zeta_1=0 \), confirming that the
  quotient geometry is non-degenerate.
  
  In the case \( u_1=1 \), the induced map \( \tilde\phi\colon \mathbb
  C^{1,1}/\mathbb T^1 \to K=K_1 \) is \( 2 \)-to-\( 1 \), branched over \(
  \partial K = W_1 \).  Each disc \( \mathbb D(r) = \{ (a,b) \in K :
  a=r+\lambda_1^{(1)} \} \) is the image of a two-sphere~\( S^2(r) \) in \(
  \mathbb C^{1,1}/\mathbb T^1 \cong \mathbb R^3 \), and this sphere in turn
  is the quotient of the three-sphere \( S^3(2r)=\{ \abs z^2+\abs w^2=2r \}
  \) in~\( \mathbb C^{1,1} \).  The map \( S^2(r)\to \mathbb D(r) \) may be
  thought of as orthogonal projection to the equatorial plane, whereas the
  map \( S^3(2r) \to S^2(r) \) is the Hopf fibration.
\end{example}

\begin{example}
  If \( n=1 \) and \( d=2 \) we are considering hypersymplectic quotients
  of \( \mathbb C^{2,2} \) by a one-dimensional Abelian group and the
  result is a four-dimensional hypersymplectic manifold~\( M \) with \( S^1
  \)-symmetry.  The map \( \beta \) defining~\( N \) is determined by \(
  u_1,u_2 \in \mathbb Z \) not both zero.  Without loss of generality we
  may take \( u_1>0 \).  Let us restrict to the non-degenerate case where
  \( u_2 \) is also non-zero.  The moment map~\( \phi \) on~\( M \) has
  image \( K=K_1\cap K_2 \), where \( K_k \) are solid cones in \( \mathbb
  R\times \mathbb C \) with vertices \( V_k =
  (\tilde\lambda^{(1)}_k,\tilde\lambda^{(c)}_k) :=
  (\lambda^{(1)}_k/u_k,\lambda^{(c)}_k/u_k) \).  The quotient \( M \)~is
  compact if and only if~\( u_1>0>u_2 \).

  Let us successively consider the conditions (F), (S) and~(D)
  of~\S\ref{sec:non-degeneracy} for these quotients.

  \begin{figure}[tbp]
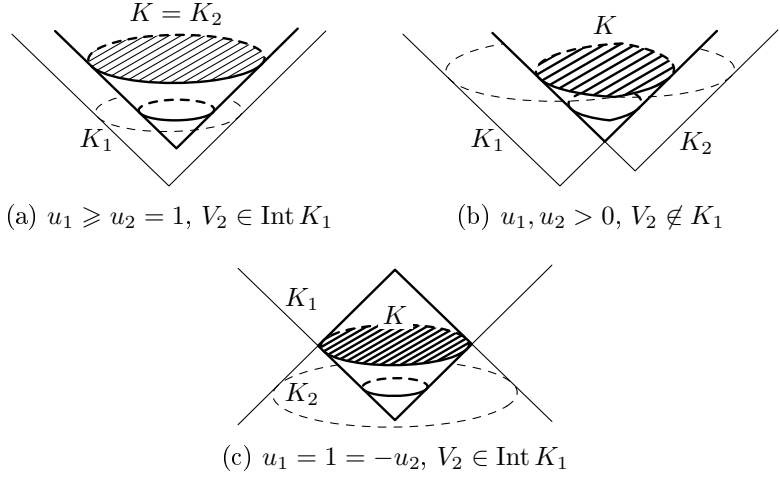

    \centering
    \subfloat[\( u_1 \geqslant u_2=1 \), \( V_2\in\Int K_1 \)]{\label{fig:d2-a}%
    \quad\includegraphics{hs-toric-pictures.2}\quad}
    \subfloat[\( u_1,u_2>0 \), \( V_2\not\in K_1 \)]{\label{fig:d2-b}%
    \quad\includegraphics{hs-toric-pictures.3}\quad}
    \subfloat[\( u_1=1=-u_2 \), \( V_2\in\Int K_1 \)]{\label{fig:d2-c}%
    \quad\includegraphics{hs-toric-pictures.4}\quad}
    \caption{Configurations giving smooth quotients from non-degenerate \(
    S^1 \)- and \( S^1\times\mathbb Z/m \)-actions on~\( \mathbb C^{2,2}
    \).}
    \label{fig:d2}
  \end{figure}

  For the freeness condition~(F) of Definition~\ref{def:free},
  Proposition~\ref{prop:stabiliser} forces \( V_1\ne V_2 \) and imposes the
  restriction that \( u_k=\pm1 \) if the vertex \( V_k \) lies in~\( K \).
  Allowable cone configurations include the three given in
  Figure~\ref{fig:d2}: \eqref{fig:d2-a}~has \( u_1>0 \), \( u_2=1 \) and \(
  N=S^1 \); \eqref{fig:d2-b},~\( u_1,u_2>0 \) and \( N=S^1\times \mathbb
  Z/\gcd(u_1,u_2) \); \eqref{fig:d2-c},~\( u_1=1 \), \( u_2=-1 \).  Other
  allowable configurations have \( V_2\in W_1\setminus V_1 \).

  Turning to the smoothness condition~(S) of Definition~\ref{def:smooth},
  assume that the \( N \)-action is free.  By
  Remark~\ref{rem:smooth-bound}, smoothness fails if \( V_2 \) lies in~\(
  W_1 \), since at~\( V_2 \) we then have \( L=\{1,2\} \), \( J=\{2\} \)
  and \( \lie n_L=\lie n \) has dimension~\( 1 \).  Thus the configurations
  of Figure~\ref{fig:d2} are the only candidates for smooth quotients.  For
  Figure~\ref{fig:d2-a}, there is no more to check as \( \lie n_L=\{0\} \)
  at all points of~\( K \).  For Figures~\ref{fig:d2-b} and~\ref{fig:d2-c}
  we need to consider points \( (a,b) \in W_1\cap W_2 \).  Here \(
  L=\{1,2\} \) and \( J=\varnothing \), so \( \lie n_{L,J} = \lie n_L =
  \lie n \cong \mathbb R \) and the map \( \Lambda_{(a,b)}\colon \lie
  n\otimes (\mathbb R\times \mathbb C)\to \mathbb C^2 \) in
  equation~\eqref{eq:Lambda} is
  \begin{equation*}
    \Lambda_{(a,b)}
    \left(
      \begin{pmatrix}
        x_1\\
        x_2
      \end{pmatrix}
    \otimes (\theta^{(1)},\theta^{(c)})
    \right)
    =
    \begin{pmatrix}
      x_1(b_1\theta^{(1)}+a_1\theta^{(c)})\\
      x_2(b_2\theta^{(1)}+a_2\theta^{(c)})
    \end{pmatrix}
    ,
  \end{equation*}
  where \( x_1u_1+x_2u_2=0 \).  This is not injective if and only if \(
  b_1/a_1 = - \theta^{(c)}/\theta^{(1)} = b_2/a_2 \) for some \(
  (\theta^{(1)},\theta^{(c)}) \).  On \( W_k \), \( a_k=\abs{b_k} \), so \(
  b_k/a_k = e^{\ii\varphi_k} \).  Looking from above, we have
  Figures~\ref{fig:discs-a} and~\ref{fig:discs-b}.
  \begin{figure}[tbp]
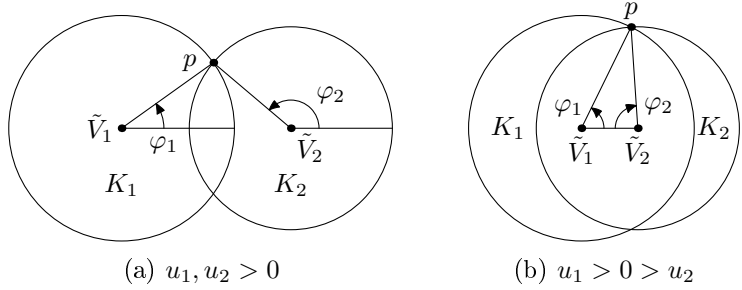

    \centering
    \subfloat[\( u_1, u_2 > 0
    \)]{\label{fig:discs-a}\includegraphics{hs-toric-pictures.5}}
    \qquad
    \subfloat[\( u_1 > 0 > u_2
    \)]{\label{fig:discs-b}\includegraphics{hs-toric-pictures.6}}
    \caption{Level sets \( a \) constant through \( K_1 \) and \( K_2 \)
    showing the angles~\( \varphi_k \) for a point \( p=(a,b)\in W_1\cap
    W_2 \).  The points \( \tilde V_k \) are the orthogonal projections to
    this plane of the vertices~\( V_k \).}
    \label{fig:discs}
  \end{figure}
  One sees that the quotient is smooth when we have \( V_k\notin W_j \) for
  \( k\ne j \).  Thus the three configurations of Figure~\ref{fig:d2} give
  smooth four-dimensional manifolds.

  Finally, we turn to the non-degeneracy of the hypersymplectic structure
  as guaranteed by condition~(D) of Definition~\ref{def:non-deg}, when the
  quotient is smooth satisfying conditions (F) and~(S).
  Remark~\ref{rem:s-zero} shows that for (D) to hold we must have \(
  W_1\cap W_2 = \varnothing \).  This only occurs when \( u_2=1 \) and \(
  V_2\in \Int K_1 \), as in Figure~\ref{fig:d2-a}.  Thus the other two
  configurations give smooth manifolds with a hypersymplectic structure
  that degenerates along some hypersurface.  Indeed in these two cases the
  fixed point set of the involution~\( \hat\sigma \)
  of~\S\ref{sec:involutions} is non-empty and
  Proposition~\ref{prop:fixed-deg} applies.

  For \( u_2=1 \) and \( V_2\in\Int K_1 \) we need to consider
  Proposition~\ref{prop:phi-non-deg} in detail.  Put \( \tilde a_k =
  a_k/u_k = a- \tilde\lambda^{(1)}_k \), \( \tilde b_k = b_k/u_k = b-
  \tilde\lambda^{(c)}_k \) and
  \begin{equation*}
    f_k(a,b) = {\tilde a_k}^2 - \abs{\smash[t]{\tilde b_k}}^2.
  \end{equation*}
  Then \( f_k \geqslant 0 \) on~\( K_k \) with equality on~\( W_k \).
  Equations~\eqref{eq:proportionality} and~\eqref{eq:lin-dep} are now
  \begin{gather}
    \zeta_1^2f_1 = \tfrac14 s^2 = \zeta_2^2f_2,\label{eq:f1-f2}\\
    \zeta_1u_1+\zeta_2u_2 = 0,\label{eq:u1-u2}
  \end{gather}
  for some \( s \in \mathbb R \).  In \( K=K_2 \), we have \( f_1>f_2 \),
  since this holds on~\( W_2 \) and \( f_1-f_2 \) is an increasing linear
  function of~\( a \).  Fixing~\( b \), we have \( \lim_{a\to\infty}
  (f_2/f_1)(a,b) = 1 \), so \( f_2/f_1 \) takes all values in \( [0,1) \)
  in~\( K \).  Now \eqref{eq:f1-f2} implies \( \zeta_1^2/\zeta_2^2 =
  f_2/f_1 \), whereas \eqref{eq:u1-u2} gives \( \zeta_1/\zeta_2 = - 1/u_1
  \), since \( u_2=1 \).  We conclude that the hypersymplectic structure on
  this quotient is non-degenerate if and only if \( u_1=1 \).
  
  Topologically, the quotient~\( M \) from \( u_1=u_2=1 \) and \(
  V_2\in\Int K_1 \), is non-compact with two connected components
  interchanged by the involution \( \hat\sigma \)
  of~\S\ref{sec:involutions}.  On each component, \( a \)~is a Morse
  function with a single critical point of index~\( 0 \).  Topologically
  and smoothly each component of~\( M \) is~\( \mathbb R^4 \).  The
  hypersymplectic structure is not in general flat.
  
  To see this, consider the case where \( \lambda_1^{(c)} = 0 =
  \lambda_2^{(c)} \), so the quotient carries the Harada-Proudfoot
  circle action of Remark~\ref{rem:circleaction}.  The fixed-point set has,
  as one component, the image of \( w_1 = 0 = w_2 \).  This is totally
  geodesic in~\( M \) with metric~\( h \) obtained from the quotient of the
  set \( \{\, \abs{z_1}^2-\abs{z_2}^2 = c_1 = \lambda_1^{(1)} +
  \lambda_2^{(2)}\, \} \) in~\( \mathbb C^2 \) by the isometric circle
  action \( (z_1,z_2)\mapsto (e^{\ii \theta}z_1,e^{-\ii \theta}z_2) \).
  For the case \( c_1=1 \), writing \( z_1 = \sqrt{1+r^2} e^{\ii\theta} \)
  and \( z_2 = r e^{\ii(\psi-\theta)} \), we get \(
  h=\frac{2r^2+1}{r^2+1}dr^2+ \frac{r^2(r^2+1)}{2r^2+1}d\psi^2 \), which
  has curvature \( -1/(2r^2+1)^3 \).
\end{example}

\subsection{Another circle action}
Let us modify the action of \S\ref{sec:diagonal}, so \( d = n+1 \) and
\begin{gather*}
  u_k = e_k, \qquad\text{for \( k=1, \dots, n \),} \\
  u_{n+1} = e_1 + \dots + e_n.
\end{gather*}
Take \( \lambda^{(i)}_k = 0 \) for \( i=1,2,3 \) and \( k=1,\dots,n \).
Put \( \lambda^{(1)}_{n+1}=-\lambda<0 \) and \( \lambda^{(c)}_{n+1}=0 \).
Then if \( a=(a_{(1)},\dots,a_{(n)}) \in \mathbb R^n \), etc., we have
\begin{gather*}
  K_k = \{\, (a,b)\in \mathbb R^n\times \mathbb C^n : a_{(k)}\geqslant
  \abs*{b_{(k)}} \,\},\qquad\text{for \( k=1,\dots,n \),}\\
  K_{n+1} = \Biggl\{\, (a,b)\in \mathbb R^n\times \mathbb C^n :
  \lambda+\sum_{k=1}^n a_{(k)} \geqslant \abs*[\bigg]{\sum_{k=1}^n b_{(k)}}
  \,\Biggr\}. 
\end{gather*}
We see immediately that \( K \) is the intersection of just \(
K_1,\dots,K_n \) and that \( W_{n+1} \) does not meet \( K \).  As \(
V_1\cap\ldots\cap V_n \) lies in \( K \) and \( \{e_1,\dots,e_n\} \) is a \(
\mathbb Z \)-basis for~\( \mathbb Z^n \), Proposition~\ref{prop:stabiliser}
implies that the action of the circle~\( N \) is free on \( \mu^{-1}(0) \).
For smoothness, note that the index sets \( L \) and \( J \) are both
subsets of \( \{1,\dots,n\} \).  However, \( \lie n_{\{1,\dots,n\}} = \{0\}
\) since the restriction of \( \beta \) to \( \mathbb
R_{\{1,\dots,n\}}=\mathbb R^n\to\mathbb R^n \) is the just the identity
map.  Proposition~\ref{prop:smoothness} implies that the hypersymplectic
quotient \( M=\mu^{-1}(0)/N \) is thus smooth.  Topologically and smoothly
it has two connected components which are copies of \( \mathbb R^{4n} \).
These components are interchanged by the involution~\( \sigma \).  

That the quotient has a non-degenerate hypersymplectic structure may be
seen as follows.  Using \eqref{eq:lin-dep}, we have \( \zeta_1 = \dots =
\zeta_n = \zeta \) and \( \zeta_{n+1} = -\zeta \).  Putting \( \tilde s =
s/{2\zeta} \) the system~\eqref{eq:proportionality} becomes
\begin{gather*}
  a_{(k)}^2 = \abs*{b_{(k)}}^2+\tilde s_k^2,\qquad\text{for \( k=1,\dots,n
  \)},\\
  (a_{(1)}+\dots+a_{(n)}+\lambda)^2=\abs*{b_{(1)}+\dots+b_{(n)}}^2
  + (\tilde s_{(1)}+\dots+\tilde s_{(n)})^2.
\end{gather*}
Using the triangle inequality, the first \( n \)-equations imply that the
left-hand side of the last equation is strictly greater than the right-hand
side for \( \lambda>0 \).  So there is no solution to the degeneracy
equations and we obtain a smooth non-degenerate hypersymplectic structure
on two copies of \( \mathbb R^{4n} \).  The computations of the previous
section show that this metric is not flat.


\begin{thebibliography}{AAAA}
  
\bibitem[AD]{AD} A. Andrada and I. Dotti, Double products and
  hypersymplectic structures on {$R^{4n}$}, January 2004, eprint
  \url{arXiv:math.DG/0401294}.
  
\bibitem[BD]{BD} R. Bielawski and A. Dancer. The geometry and topology of
  toric hyperk\"ahler manifolds. Comm. Anal. Geom.  \textbf{8} (2000)
  727--759.

\bibitem[De]{De} T. Delzant. Hamiltoniens p\'eriodiques et images convexe
  de l'application moment. Bull. Soc. Math. France \textbf{116} (1988)
  315--339.

\bibitem[FPPW]{FPPW} A. Fino, H. Pedersen, Y.~S. Poon and M. Weye S\o
  rensen.  Neutral Calabi-Yau structures on Kodaira manifolds. IMADA
  preprint PP-2002-13, Odense 2002,
  \url{http://bib.mathematics.dk/preprint.php?id=IMADA-PP-2002-13},
  Commun. Math. Phys. (to appear).

\bibitem[Ga]{Ga} K. Galicki. A generalisation of the momentum mapping
  construction for quaternionic K\"ahler manifolds.  Commun. Math. Phys.
  \textbf{108} (1987) 117--138.

\bibitem[GL]{GL} K. Galicki and H. B. Lawson. Quaternionic reduction and
  quaternionic orbifolds. Math. Annalen \textbf{282} (1988) 1--21.

\bibitem[Gu]{Gu} V. Guillemin. K\"ahler structures on toric varieties. J.
  Diff. Geom. \textbf{40} (1994) 285--309.

\bibitem[HP]{HP} M. Harada and N. Proudfoot. Properties of the residual
  circle action on a toric hyperk\"ahler variety. eprint
  \url{arXiv:math.DG/0207012}.

\bibitem[HKLR]{HKLR} N.~J. Hitchin, A. Karlhede, U. Lindstr\"om and M.
  Ro\v cek.  Hyperk\"ahler metrics and supersymmetry.  Commun. Math. Phys.
  \textbf{108} (1987) 535--589.

\bibitem[H]{H} N.~J. Hitchin. Hypersymplectic quotients. Acta Acad. Sci.
  Tauriensis, supplemento al numero \textbf{124} (1990) 169--180.

\bibitem[Hu]{Hu}
  C.~M. Hull.  Actions for $(2,1)$ sigma models and strings. Nuclear Phys.
  B \textbf{509} (1998), no.~1-2, 252--272.

\bibitem[Ka]{Ka} H.~Kamada.  Neutral hyper-K\"ahler structures on primary
  Kodaira surfaces.  Tsukuba J. Math. \textbf{23} (2) (1999) 321--332.

\end{thebibliography}
\end{document}